\documentclass[11pt]{article}
\usepackage[margin=1in]{geometry}
\usepackage{graphicx}
\usepackage{booktabs}
\usepackage{longtable}
\usepackage{array}
\usepackage{float}
\usepackage{caption}
\usepackage{enumitem}
\usepackage[hidelinks]{hyperref}
\usepackage[super,sort&compress]{natbib}
\usepackage{xcolor}
\usepackage{setspace}
\usepackage{titlesec}
\titleformat{\section}{\large\bfseries}{\thesection.}{0.6em}{}
\titleformat{\subsection}{\normalsize\bfseries}{\thesubsection}{0.6em}{}
\setlist[itemize]{leftmargin=1.25em}
\usepackage{varioref}
\usepackage{hyperref}
\hypersetup{
	colorlinks=true,
	linkcolor=violet,
	citecolor=cyan,
	urlcolor=teal
}
\usepackage{authblk}
\usepackage{titlesec}

\newtheorem{Rem}{Remark}

\graphicspath{{Figure/}}
\date{}
\begin{document}
	%\title{Measles Resurgence in Bangladesh, 2026: A Rapid Review and Situational Analysis for Urgent Public Health Response}
	\title{Measles Resurgence in Bangladesh, 2026: A Situational Analysis for Urgent Public Health Response}
	
	\author[1]{\small Md. Kamrujjaman}%\thanks{Corresponding author   Email: kamrujjaman@du.ac.bd} }
	\author[1, 2]{\small Faizunnesa Khondaker}%\thanks{Email: faizunnesa@math.jnu.ac.bd}}
	\author[1]{\small Alvi Ahmed Sarker}%\thanks{Email: a.a.sarker09@gmail.com} }
	\author[1,3]{\small  Nuzhat Nuari Khan Rivu}%\thanks{Email: nuzhatnuarikhan-rivu@uiowa.edu} }

	%	\author[2]{\small  Md. Shahidul Islam\thanks{Email: mshahid@du.ac.bd}}

	\affil[1]{\footnotesize Department of Mathematics, University of Dhaka, Dhaka 1000, Bangladesh}
	\affil[2]{\footnotesize Department of Mathematics, Jagannath University, Dhaka 1100, Bangladesh}
	\affil[3]{\footnotesize Department of Mathematics, University of Iowa, Iowa City, IA, 52242, USA}
	\affil[ ]{\textbf{Corresponding author:} Md. Kamrujjaman, \href{mailto:kamrujjaman@du.ac.bd}{kamrujjaman@du.ac.bd}}
	\maketitle

\begin{abstract}
	\noindent\textbf{Background} Measles has returned with force in the post-pandemic period as routine immunisation recovery has stalled well below the two-dose coverage required to interrupt transmission. %The World Health Organization reported 10.3 million cases globally in 2023, and the South-East Asia Region recorded a 42\% increase in reported cases in 2025 compared with the previous five-year average. 
	Bangladesh, which had moved close to measles--rubella elimination, entered 2026 with widening coverage gaps, depleted vaccine stocks, and mounting concern about missed children. We undertook a rapid review and situation analysis to clarify the scale, concentration, and programme implications of the outbreak while response decisions were still unfolding.
	
	\medskip
	\noindent\textbf{Methods} We conducted a rapid mixed-evidence review from 1 to 15 April 2026 across WHO, UNICEF, Directorate General of Health Services (DGHS) bulletins, PubMed/MEDLINE, ReliefWeb, WHO South-East Asia Regional Office (SEARO) updates, and Bangla-language media. We screened 46 records and included 15 primary outbreak sources and four contextual peer-reviewed references. The situation analysis used aggregated, publicly available surveillance and programme data. %We summarised temporal, geographic, demographic, and vaccination-status patterns and quantified division-level concentration with Pareto analysis and the Herfindahl--Hirschman Index (HHI). %Suspected and confirmed case definitions followed national and WHO guidance.
	
	\medskip
	\noindent\textbf{Findings} By 15 April 2026, Bangladesh had recorded 19,161 suspected measles cases, 2,973 confirmed cases, 166 suspected deaths, and 32 confirmed deaths across 58 districts since 15 March 2026. Burden was concentrated rather than diffuse: the two highest-burden divisions accounted for 56.5\% of reported cases in the division-level reporting slice analysed, and the HHI was 0.217. Children younger than five years represented 81\% of reported cases, including 34\% in infants younger than nine months. Reported vaccination status showed 72\% zero-dose and 16\% partially vaccinated cases. Programme indicators pointed to declining protection, with valid measles-rubella first-dose coverage falling from 88.6\% in 2019 to 86\% in 2024 and second-dose coverage from 89\% to 80.7\%, leaving an estimated 20 million children vulnerable on the second-dose gap measure.
	
	\medskip
	\noindent\textbf{Interpretation} We found that the 2026 resurgence reflected accumulated immunity gaps rather than biological vaccine failure. The outbreak exposed the consequences of subnational inequity, zero-dose accumulation, programme disruption, and interrupted vitamin A delivery. Immediate priorities are upazila-level targeting of the emergency measles--rubella campaign, restoration of vitamin A supplementation, paediatric surge capacity in high-burden hospitals, and integration of real-time surveillance into outbreak operations.
	
	%\medskip
	%\noindent\textbf{Funding} [This research received no specific funding from any funding agency in the public, commercial, or not-for-profit sectors. The author M. Kamrujjaman acknowledged to the University Grants Commission (UGC), Bangladesh.]
\end{abstract}

\clearpage
\begin{Rem}
While we are preparing this report for submission, the situation in Bangladesh remains highly concerning. Recent reports (up to April 22, 2026) indicate that, on average, approximately six children are dying each day, see \url{http://dghs.portal.gov.bd/pages/press-releases} %, and the number of infections is doubling approximately every five days .
\end{Rem}

%\clearpage
\medskip
\noindent\textbf{Highlights}
\begin{itemize}
	
	\item Measles has resurged globally in the post-pandemic period, with 10.3 million cases reported in 2023 and a 42\% increase in South-East Asia in 2025.
	
	\item Bangladesh entered 2026 with widening immunisation gaps, declining vaccine coverage, and growing numbers of zero-dose and under-immunised children.
	
	\item By 15 April 2026, 19{,}161 suspected and 2{,}973 confirmed measles cases, along with 166 suspected and 32 confirmed deaths, were reported across 58 districts.
	
	\item The outbreak burden was highly concentrated, with the two most affected divisions accounting for 56.5\% of cases (HHI = 0.217).
	
	\item Children under five accounted for 81\% of cases, including 34\% among infants younger than nine months.
	
	\item A large majority of cases occurred among unvaccinated (72\%) or partially vaccinated (16\%) individuals.
	
	\item Measles--rubella vaccine coverage declined between 2019 and 2024, leaving an estimated 20 million children vulnerable due to second-dose gaps.
	
	\item The resurgence reflects accumulated immunity gaps rather than vaccine failure, driven by subnational inequities and programme disruptions.
	
	\item Immediate priorities include targeted emergency vaccination campaigns, restoration of vitamin A supplementation, strengthening paediatric care capacity, and integrating real-time surveillance into outbreak response.
\end{itemize}
\clearpage
\section{Introduction}
Measles returned to Bangladesh in the spring of 2026 with a brutality that no country approaching elimination should accept. Wards filled with coughing infants, isolation beds overflowed, and a disease long described as preventable again became a cause of avoidable child death. This happened in a world already on warning. The Immunization Agenda 2030 mid-point review made clear that global recovery after COVID-19 remained incomplete, while WHO and UNICEF reported that first- and second-dose measles coverage still sat below the 95\% threshold needed to interrupt sustained transmission. WHO estimated 10.3 million cases globally in 2023, and WHO later reported about 95,000 measles deaths in 2024. In South-East Asia, reported measles cases in 2025 rose 42\% above the previous five-year average. The global context matters because Bangladesh's outbreak did not emerge in isolation; it unfolded as part of a broader return of measles wherever immunity gaps were allowed to accumulate \citep{WHO_IA2030_2025,WHO_measlesdeaths_2025,WUENIC_2025,WHO_SEARO_EpiBull_2026}.

Bangladesh had earned a different reputation. Since nationwide expansion of the Expanded Programme on Immunisation in 1985, the country built one of the region's most respected vaccine delivery systems. Supplementary measles--rubella campaigns, systematic outreach, and introduction of the routine second dose in 2012 pushed the country close to elimination and sharply reduced disease burden. Yet Bangladesh's success was never simply a matter of national averages. It depended on whether the programme could keep reaching children in urban informal settlements, mobile communities, hard-to-reach districts, and refugee settings where missed doses accumulate quietly until transmission exposes them. That fragility had been documented before 2026, including in Rohingya and host-community settings in Cox's Bazar \citep{Khanal_MMWR_2017,Uddin_BMCInfDis_2016,WHO_Leaving_2020,WUENIC_BGD_2024,feldstein_rohingya_sero_2020}.

The paradox of 2026, then, was not that measles returned. It was that it returned after years of apparent control and did so with extraordinary speed. The available evidence points to converging failures rather than a single trigger: routine measles--rubella coverage weakened through 2024 and 2025; no national supplementary campaign had been conducted since 2020; stocks were reported as depleted; and political transition disrupted programme continuity when catch-up efforts should have accelerated. Early signals in Cox's Bazar and Dhaka's crowded settlements were followed by rapid spread across much of the country. In this setting, measles did exactly what it always does when susceptibility builds up: it found children who had never been vaccinated, children who had started but not completed the schedule, and infants too young to be protected by the routine first dose \citep{LancetMicrobe_GlobalResurgence_2026,LancetSEA_Elimination_2025,Prothomalo_March2026,BMJ_Bangladesh_2026,UNICEF_SitRep1_2026,BusinessStandard_VaxShortage_2026,Guardian_Bangladesh_2026}.

We therefore undertook a rapid review and data-driven situation analysis of the 2026 measles resurgence in Bangladesh. Our aim was practical rather than retrospective. We wanted to synthesise the best available evidence while the outbreak was still evolving, show where the burden was concentrating, identify which children were being missed, and translate that evidence into immediate actions for Bangladeshi health authorities, programme managers, and regional partners facing the same post-pandemic vulnerabilities across South Asia \citep{LancetSEA_Elimination_2025,WHO_IA2030_2025}.

\section{Measles Epidemiology and Outbreak Response in Bangladesh}
\subsection{Global Measles Situation Across WHO Regions}
\begin{figure}[H]
\centering
\includegraphics[width=0.99\textwidth]{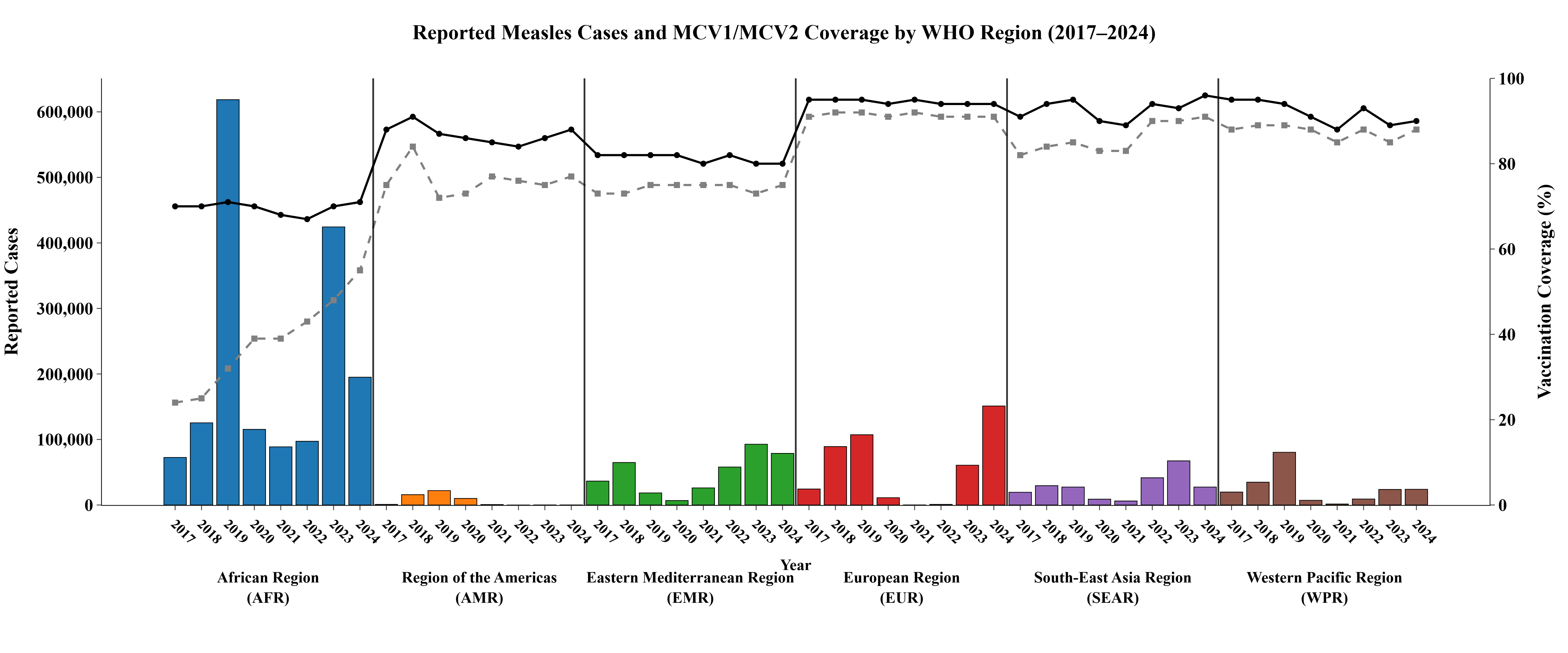}
\caption{Reported measles cases and first-dose and second-dose measles-containing vaccine coverage by WHO region, 2017--2024. The figure highlights marked regional differences in reported measles burden despite generally high routine coverage in several regions. Source: WHO Immunization Data Portal provisional measles and rubella data.}\label{fig:who_regions}
\end{figure}

The measles situation is not a localized fluctuation but a sustained global resurgence (Figure~\ref{fig:who_regions}). The {World Health Organization} (WHO) and the {U.S. Centers for Disease Control and Prevention} (CDC) estimated approximately 10.3 million measles cases and 107{,}500 deaths worldwide in 2023. More recent WHO reports indicate an estimated 95{,}000 deaths in 2024. 
Global vaccination coverage improved only modestly in 2024, reaching 84\% for the first dose and 76\% for the second dose of the measles vaccine. These levels remain substantially below the 95\% two-dose threshold required to interrupt transmission \cite{who_factsheet2025}.

\begin{figure}[H]
\centering
\includegraphics[width=0.99\textwidth]{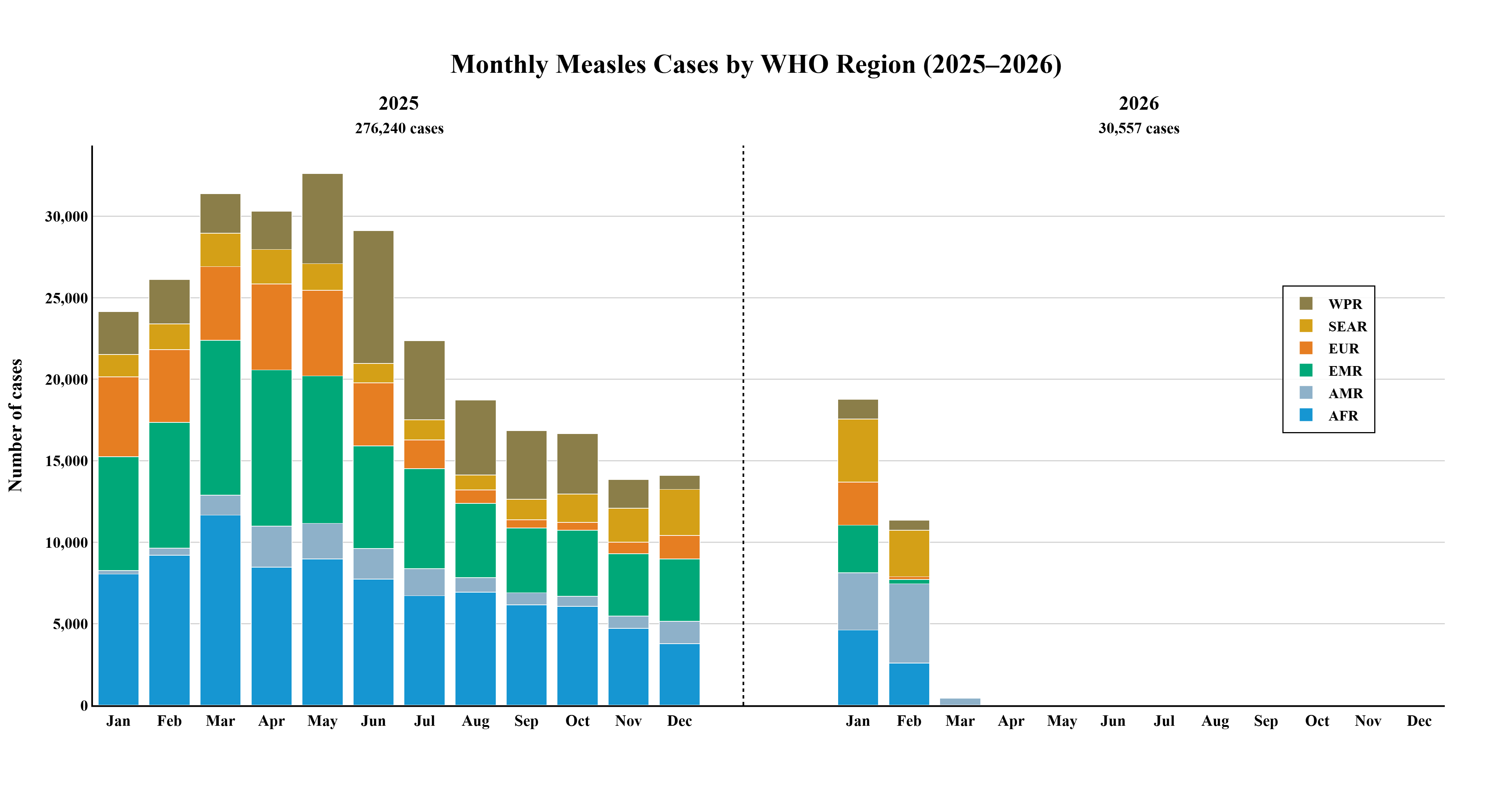}
\caption{Monthly reported measles cases by WHO region in 2025 and 2026. The early-2026 pattern shows increasing reported burden in the Americas and South-East Asia, with lower totals than the 2025 peak but persistent multi-regional activity.}\label{fig:who_months}
\end{figure}

In April 2025, WHO, UNICEF, and Gavi reported that 61 countries experienced large or disruptive measles outbreaks within the preceding 12 months. Furthermore, a CDC global update from March 2026 indicates that outbreaks are now occurring in every region worldwide \cite{surge}.
Regional trends in early 2026 show mixed patterns. During January and February 2026, reported cases increased in the Americas (AMR) and South-East Asia (SEAR) compared to the same period in 2025, while declines were observed in Africa (AFR), the Eastern Mediterranean (EMR), Europe (EUR), and the Western Pacific (WPR). Despite these regional variations, the global total in early 2026 remains lower than in 2025 (Figure~\ref{fig:who_months}).

Figure~\ref{fig:who_regions} illustrates variations in reported disease burden and routine vaccination coverage across WHO regions from 2017 to 2024. Figure~\ref{fig:who_months} highlights continued volatility in early 2026, with increases in the Americas and South-East Asia despite declines in other regions \citep{who_factsheet2025,surge,provisional,WHO_SEARO_EpiBull_2026}.

\subsection{Current Measles Situation in Bangladesh}
\begin{figure}[H]
\centering
\includegraphics[width=0.95\textwidth]{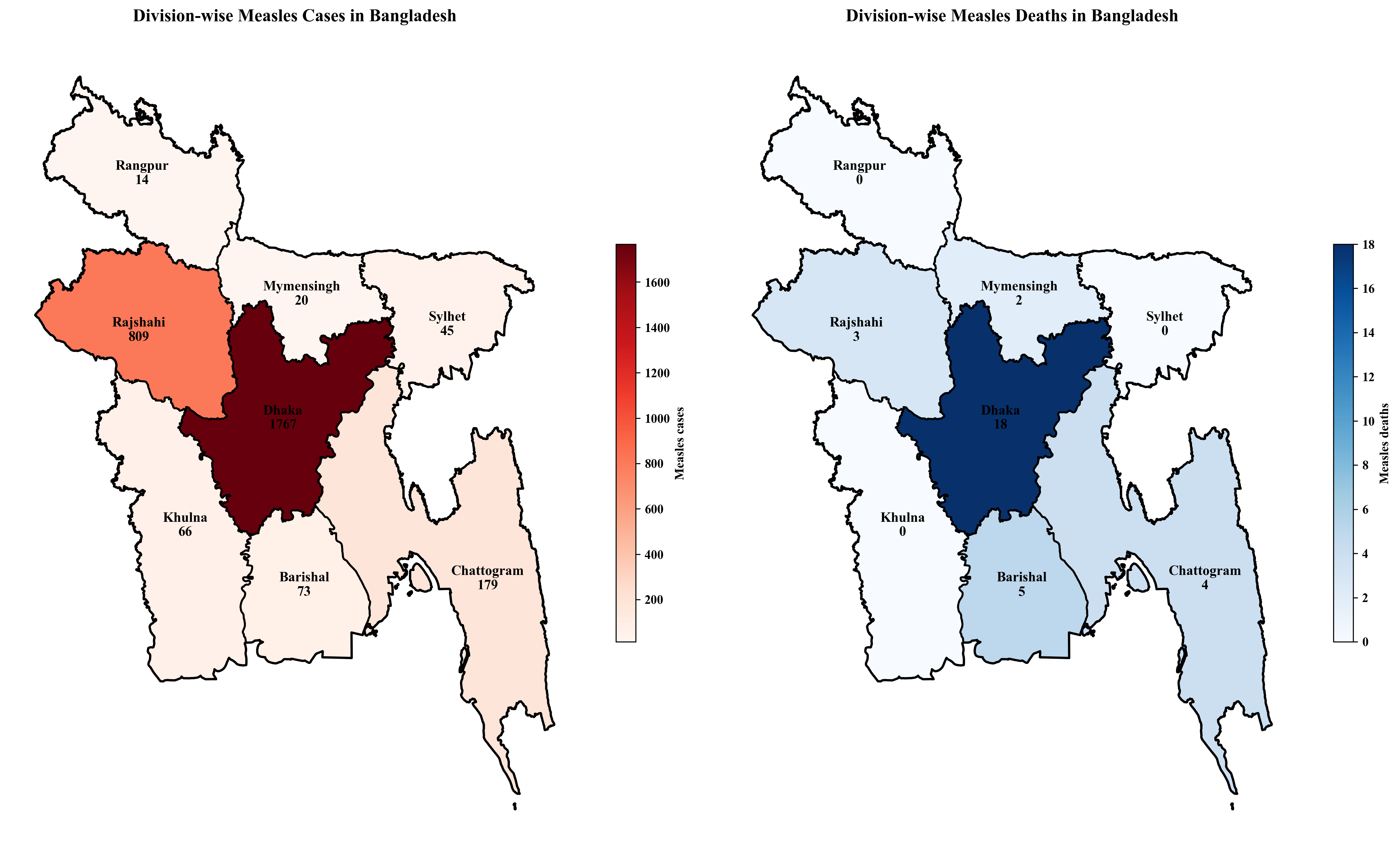}
\caption{Division-wise geographic distribution of cumulative reported measles cases and deaths in Bangladesh from 15 March to 15 April 2026. Numbers on each map indicate the total counts recorded in each division during the reporting period.}\label{fig:map_cases_deaths}
\end{figure}

\begin{figure}[H]
\centering
\includegraphics[width=0.95\textwidth]{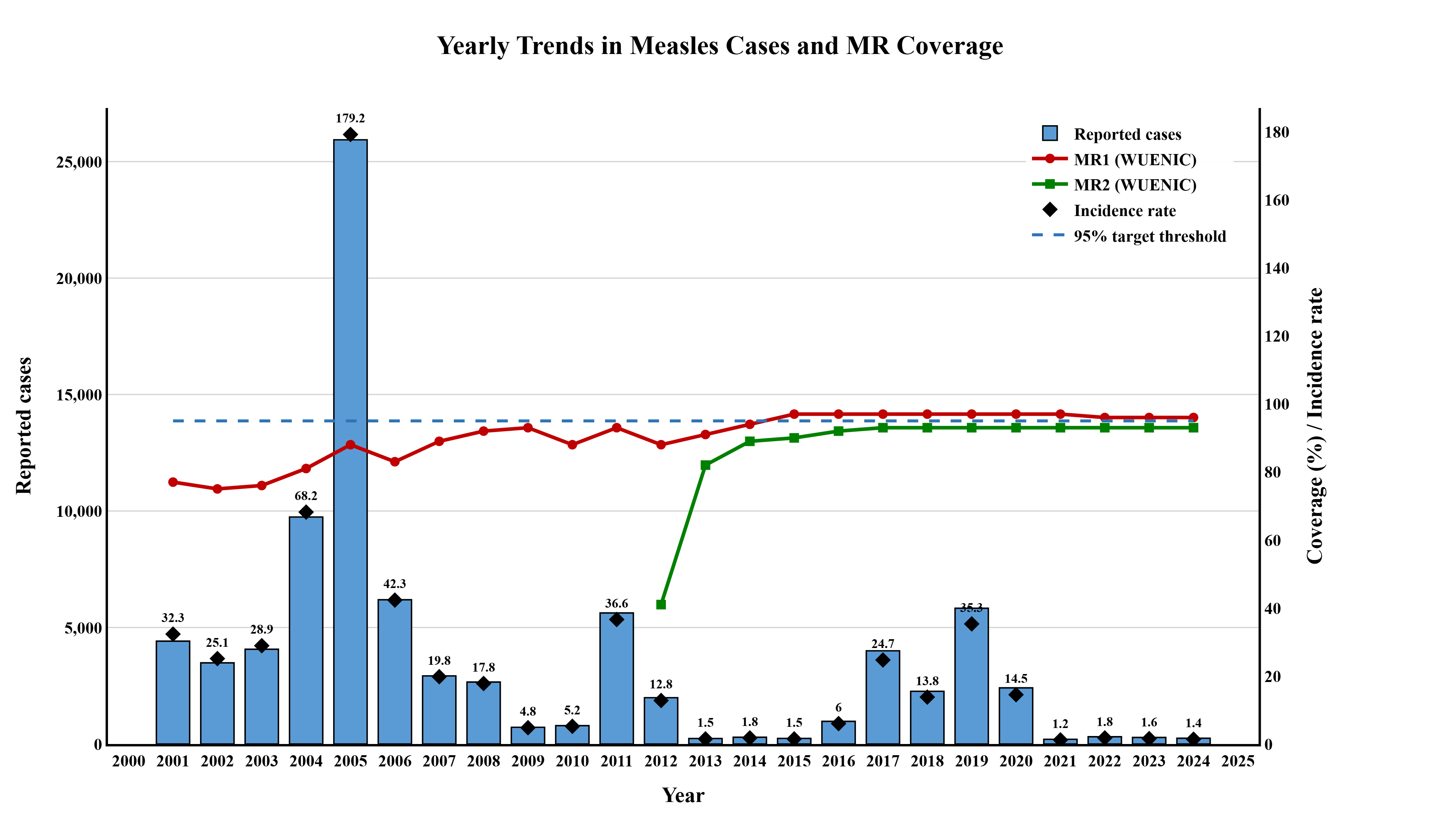}
\caption{Yearly trends in reported measles cases, incidence rate, and first-dose and second-dose measles-rubella coverage in Bangladesh. The figure shows that long-term gains in coverage broadly coincided with lower disease burden but did not prevent renewed outbreak risk when coverage and equity weakened.}\label{fig:trend_cases_coverage}
\end{figure}

Bangladesh is experiencing a fast-growing measles outbreak that intensified from mid-March 2026, with rapid geographic spread and high pediatric mortality reported through national surveillance and partner situation updates \cite{cdn}.
As of 15 April 2026, 19,161 suspected and 2,973 laboratory-confirmed measles cases were reported since 15 March 2026, alongside 166 suspected measles-associated deaths and 32 laboratory-confirmed measles deaths (Figure~\ref{fig:map_cases_deaths})  \cite{dghs}. Bangladesh's long-run national story is one of strong immunization performance with increasingly visible pockets of failure. WHO-UNICEF national estimates for 2019--2023 keep Bangladesh at 97$\%$ MCV1 and 93$\%$ MCV2, and the WHO dashboard summary for the 2024 revision shows 96$\%$ MCV1 and 93$\%$ MCV2 (Figure~\ref{fig:trend_cases_coverage}). Those national figures, however, mask marked subnational inequity. The Coverage Evaluation Survey 2023 found valid MR1 coverage about 86$\%$ nationally and valid MR2 coverage about 80.7$\%$ nationally, while UNICEF, WHO and Gavi warned in April 2025 that Bangladesh had only 81.6$\%$ fully immunized children, with 70,000 zero--dose and 400,000 under--immunized children overall and substantially worse urban performance than rural \cite{WUENIC_BGD_2024}.

The outbreak's risk profile is strongly shaped by immunity gaps: UNICEF reports that MR1 coverage declined to 86$\%$ (from 88.6$\%$ in 2019) and MR2 declined to 80.7$\%$ (from 89$\%$ in 2019)(Figure~\ref{fig:trend_cases_coverage}), leaving about 10 million children susceptible (MR1 gap) and about 20 million susceptible (MR2 gap). These coverage levels are below the 95$\%$ two-dose coverage typically required to prevent sustained measles transmission \cite{UNICEF_FAQ_2026}.
The geographic pattern is now broad rather than localized. WHO's 15 April update described 58 affected districts and highlighted the highest cumulative caseloads in Dhaka, Rajshahi, and Chattogram Divisions. The same WHO update noted eight laboratory--confirmed cases in Cox's Bazar, including cases in Rohingya refugee camps, which matters because crowding and mobile populations can amplify transmission \cite{who_searo2026}.

\begin{figure}[H]
\centering
\includegraphics[width=0.95\textwidth]{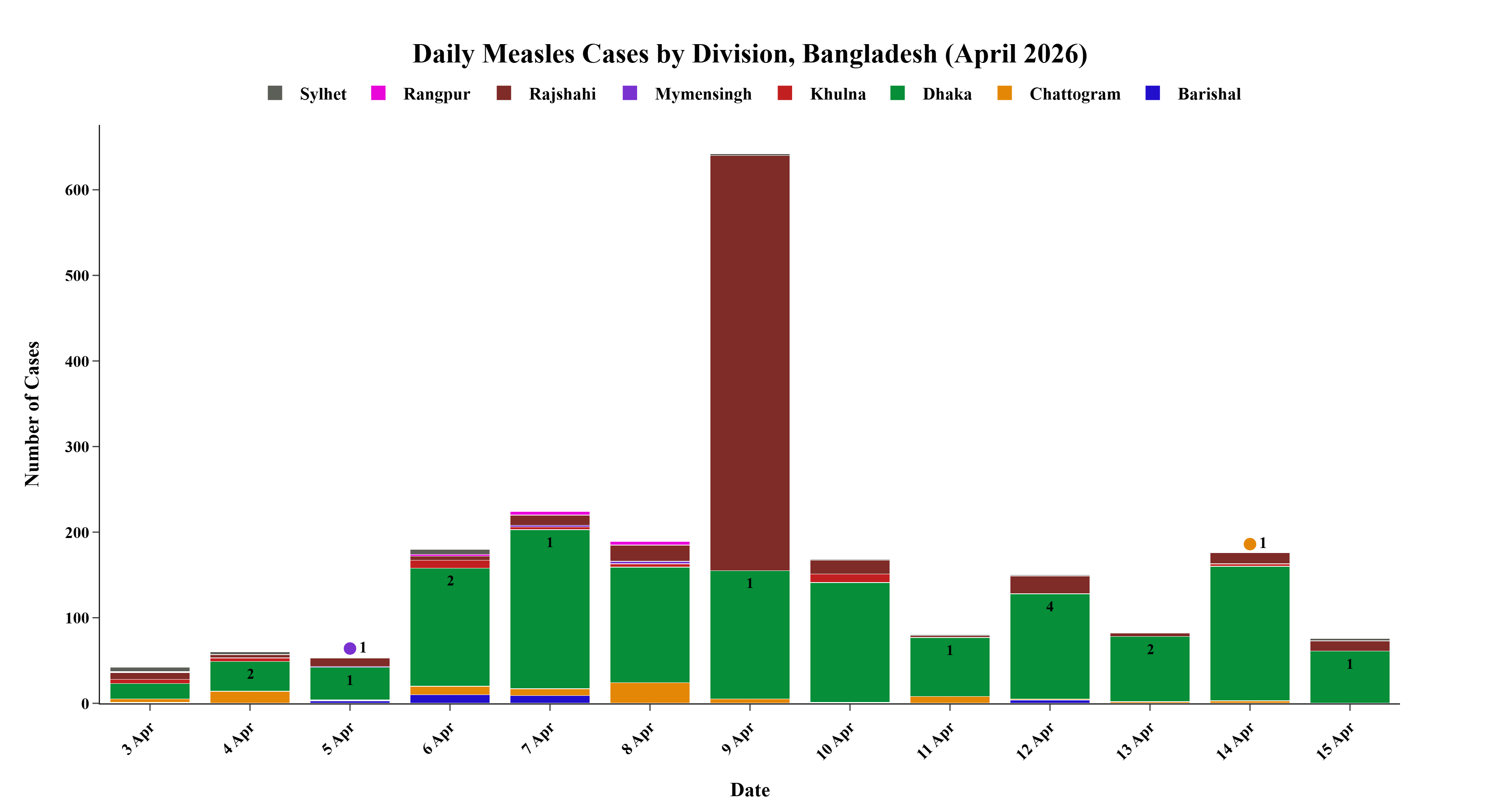}
\caption{Daily confirmed measles cases and confirmed deaths by division in Bangladesh in April 2026. Stacked bars represent confirmed cases, and labels above bars indicate confirmed deaths recorded on that day.}\label{fig:daily_confirmed_division}
\end{figure}

The outbreak is overwhelmingly pediatric. UNICEF's Bangladesh measles situation report stated that 81$\%$ of cases were in children under five, including 34$\%$ in infants under nine months. UNICEF reported that at the Infectious Diseases Hospital (IDH) in Dhaka there were 685 suspected cases, including 573 admissions in March 2026, compared with 69 in the same period of 2025; in Cox's Bazar District hospital, an eight-bed isolation unit was managing 30--40 children at a time \cite{humanitarian}.  Figure~\ref{fig:daily_confirmed_division} shows the day-to-day distribution of confirmed cases and deaths by division in April 2026. The data point in one direction: this was a fast-moving national paediatric emergency rather than a routine seasonal fluctuation \citep{who_searo2026,DGHS_Bulletins,UNICEF_SitRep1_2026,humanitarian,WUENIC_BGD_2024}.

	\subsection{Drivers of the Outbreak and Emergency Response}
The outbreak is being driven largely by immunization gaps, with many cases occurring among infants below vaccination age, along with zero-dose and under-immunized children. Although Bangladesh has historically maintained high routine immunization coverage, disruptions in MR1 and MR2 services during 2024 and 2025, despite repeated UNICEF appeals, have allowed immunity gaps to accumulate. In response, the Government of Bangladesh, with support from UNICEF, WHO, and GAVI, launched an emergency measles--rubella campaign on 5 April 2026 targeting over 1.2 million children in 30 hotspot upazilas across 18 districts; 75,442 children, or 13$\%$ of the target, were vaccinated on the first day, with expansion planned to four City Corporations from 12 April and nationwide from 3 May.
%https://www.unicef.org/media/179846/file/Bangladesh-Humanitarian-Situation-Report-No.1%28Measles-Outbreak%29-8-April-2026.pdf.pdf
The growing number of cases among infants under nine months indicates sustained community transmission and led the Government to lower the age for the first dose to six months in the emergency rollout. At the same time, disruption of Vitamin A supplementation in 2025, when only one of two rounds was completed, together with high levels of malnutrition, is increasing the risks of severe disease, complications, and death. To reduce these risks, the Ministry of Health and Family Welfare has priorotized available Vitamin A stocks for children with measles, while additional procurement and distribution are recommended to protect an estimated 20 million children under five \citep{WHO_BGD_Campaign_2026,UNICEF_SitRep1_2026,humanitarian,IFRC_FieldReport_2026}..

%\subsection{Emergency response before formal analysis}
%The emergency response reflected the scale of the immunity gap. On 5 April 2026, the Government of Bangladesh, with support from UNICEF, WHO, and Gavi, launched an emergency measles--rubella campaign targeting more than 1.2 million children in 30 hotspot upazilas across 18 districts. Officials vaccinated 75,442 children on the first day, or 13\% of the target population, and then planned expansion to city corporations from 12 April and nationwide scale-up from 3 May. The rising burden among infants younger than nine months led the government to lower the age for the first outbreak dose to six months. At the same time, health facilities were already under strain. UNICEF reported 685 suspected cases at the Infectious Diseases Hospital in Dhaka, including 573 admissions in March 2026 compared with 69 in March 2025, and an eight-bed isolation unit in Cox's Bazar was managing 30 to 40 children at a time. These operational details matter because they show that Bangladesh was attempting to close transmission gaps and mortality gaps at the same time \citep{WHO_BGD_Campaign_2026,UNICEF_SitRep1_2026,humanitarian,IFRC_FieldReport_2026}.

\section{Methods}
\subsection{Study design and rationale}
We conducted a rapid mixed-evidence review with a concurrent data-driven situation analysis between 11 and 15 April 2026. We chose this design because Bangladesh was in the middle of an active outbreak, public-health decisions could not wait for formal outbreak investigations, and the most useful evidence was already available through aggregated surveillance bulletins, partner situation reports, and a small but relevant contextual literature. Our objective was to synthesise what could be known quickly and responsibly from secondary data, not to generate individual-level causal estimates during an emergency \citep{UNICEF_SitRep1_2026,DGHS_Bulletins,WHO_SEARO_EpiBull_2026}.

\subsection{Search strategy and source selection}
We searched WHO resources, the WHO Global Immunization Data Portal, UNICEF Bangladesh situation reports, DGHS surveillance bulletins and press releases, PubMed/MEDLINE, ReliefWeb, WHO SEARO Epidemiological Bulletins, and Bangla-language media sources including Prothom Alo. Core search terms were ``measles Bangladesh 2026'', ``measles outbreak Bangladesh'', ``MR vaccine coverage Bangladesh'', and ``zero-dose children Bangladesh''. We restricted outbreak-specific sources to material published or updated between January 2020 and 15 April 2026. We screened 46 records, reviewed full texts or authoritative summaries, and included 15 primary outbreak sources plus four contextual peer-reviewed references. We included sources that focused on Bangladesh and contributed evidence on measles or measles--rubella vaccination, outbreak response, serology, vaccine effectiveness, or immunisation coverage. We excluded duplicate reports, items without identifiable institutional authorship, and non-English or non-Bangla sources that could not be translated reliably \citep{UNICEF_SitRep1_2026,WHO_SEARO_EpiBull_2026,WUENIC_BGD_2024,reliefweb_bangladesh_measles_flash_2026}.

\subsection{Data sources for the situation analysis}
The situation analysis relied on aggregated, publicly available datasets and reports. These included national daily surveillance updates from DGHS; UNICEF Bangladesh Measles Outbreak Situation Report No. 1 and Humanitarian Situation Report No. 1; the WHO SEARO epidemiological bulletin; WHO--UNICEF national immunisation coverage estimates, including the Bangladesh 2023 revision; the Bangladesh Coverage Evaluation Survey 2023 values reported in partner documents; and the ReliefWeb ECHO Daily Flash of 7 April 2026. We also used contemporaneous reporting from credible Bangla and English-language media when those reports relayed named officials, hospital observations, or figures attributable to DGHS, WHO, or UNICEF. No individual-level patient records were accessed \citep{DGHS_Bulletins,UNICEF_SitRep1_2026,humanitarian,WHO_SEARO_EpiBull_2026,WUENIC_BGD_2024,reliefweb_bangladesh_measles_flash_2026}.

\subsection{Analytical approach and case definitions}
Our analysis was descriptive. We examined temporal patterns through epidemic curves and a historical comparison using a three-year moving average; geographic patterns through division-level mapping and reporting slices; and demographic patterns through age and vaccination-status distributions. To describe geographic concentration, we ranked divisions by reported case share, calculated the cumulative share of burden, and summarised concentration with the Herfindahl--Hirschman Index, where 0 represents perfect equality and 1 full concentration. We did not conduct statistical hypothesis testing or claim causal inference. We classified a suspected measles case as fever with maculopapular rash and at least one of cough, coryza, or conjunctivitis, consistent with national and WHO guidance; confirmed cases were those with laboratory confirmation by IgM serology or PCR. Public bulletins did not consistently report probable cases separately, so we did not analyse that category. This study used only aggregated, publicly available secondary data; no ethical review or individual consent was required under standard institutional policy for secondary analyses of public surveillance information \citep{WHO_measlesdeaths_2025,DGHS_Bulletins,UNICEF_SitRep1_2026}.

\section{Results}
\subsection{Outbreak magnitude and temporal trend}
\begin{figure}[H]
\centering
\includegraphics[width=0.9\textwidth]{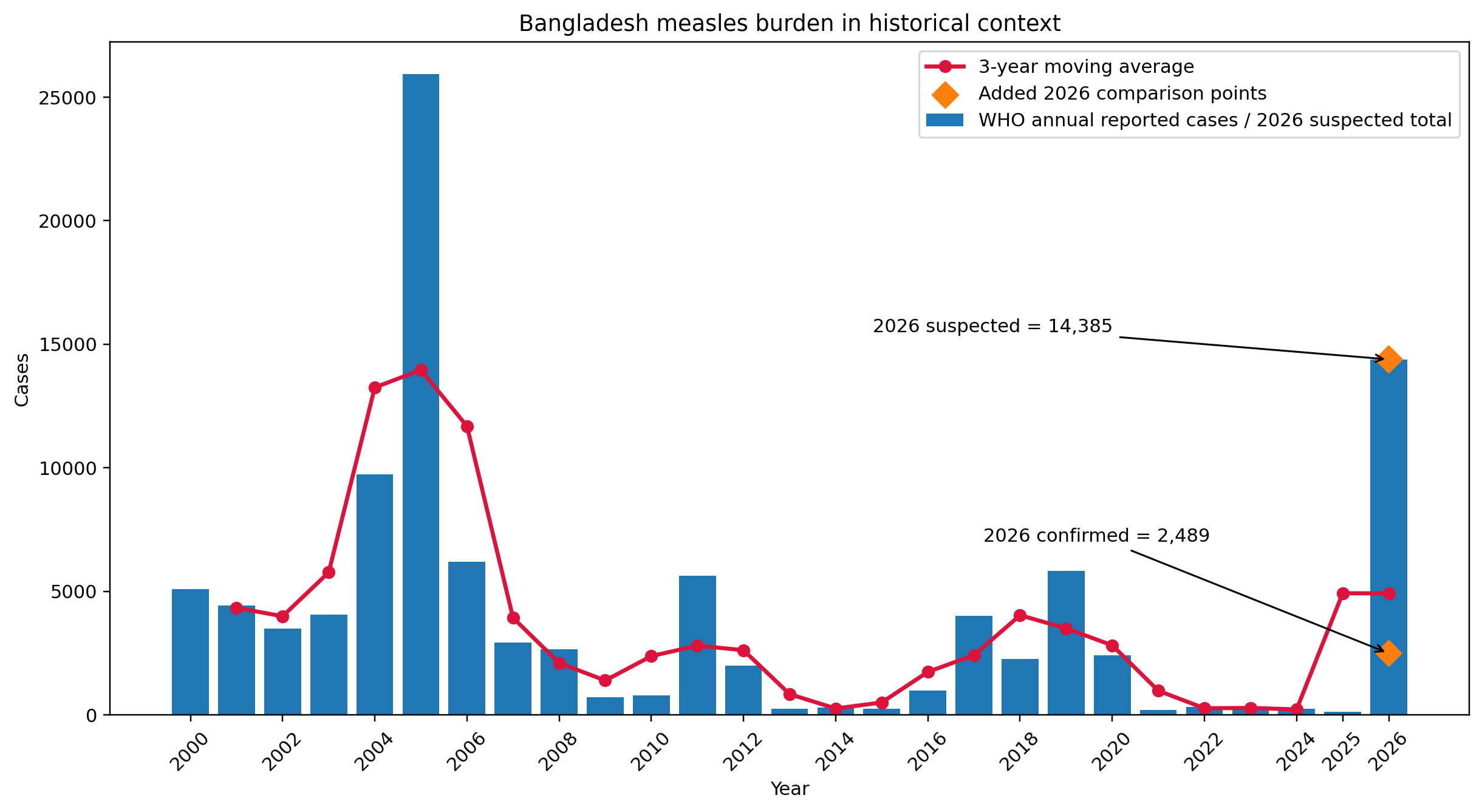}
\caption{Bangladesh measles burden in historical context. Annual WHO-reported measles cases are shown with a three-year moving average. The 2026 suspected and confirmed values are plotted separately as partial-year outbreak comparison points.}\label{fig:historical_context}
\end{figure}

Bangladesh's 2026 outbreak was historically exceptional by the middle of April. By 15 April 2026, surveillance sources had recorded 19,161 suspected cases, 2,973 laboratory-confirmed cases, 166 suspected deaths, and 32 confirmed deaths across 58 districts since 15 March 2026. Figure~\ref{fig:historical_context} shows how sharply these totals diverged from the country's recent low-burden baseline. The gap between the 2026 outbreak totals and the small annual totals reported in the years immediately preceding the outbreak was too large to dismiss as routine variation. This was a national resurgence \citep{who_searo2026,DGHS_Bulletins}.

\subsection{Geographic distribution and concentration}
\begin{figure}[H]
\centering
\includegraphics[width=0.9\textwidth]{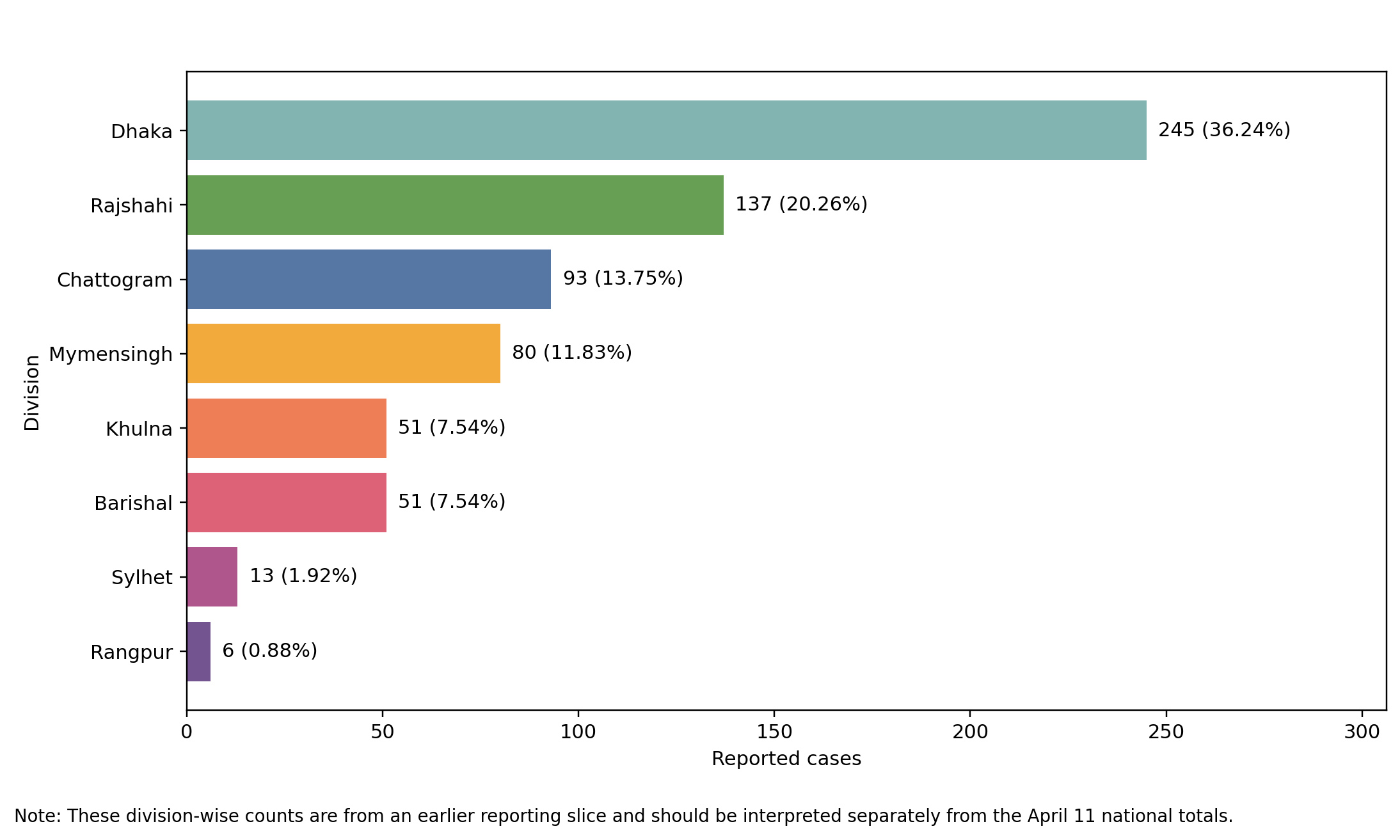}
\caption{Division-wise measles reporting in Bangladesh from the 31 March 2026 reporting slice. Dhaka accounted for the largest share of reported cases, followed by Rajshahi and Chattogram.}\label{fig:division_slice}
\end{figure}

\begin{figure}[H]
\centering
\includegraphics[width=0.88\textwidth]{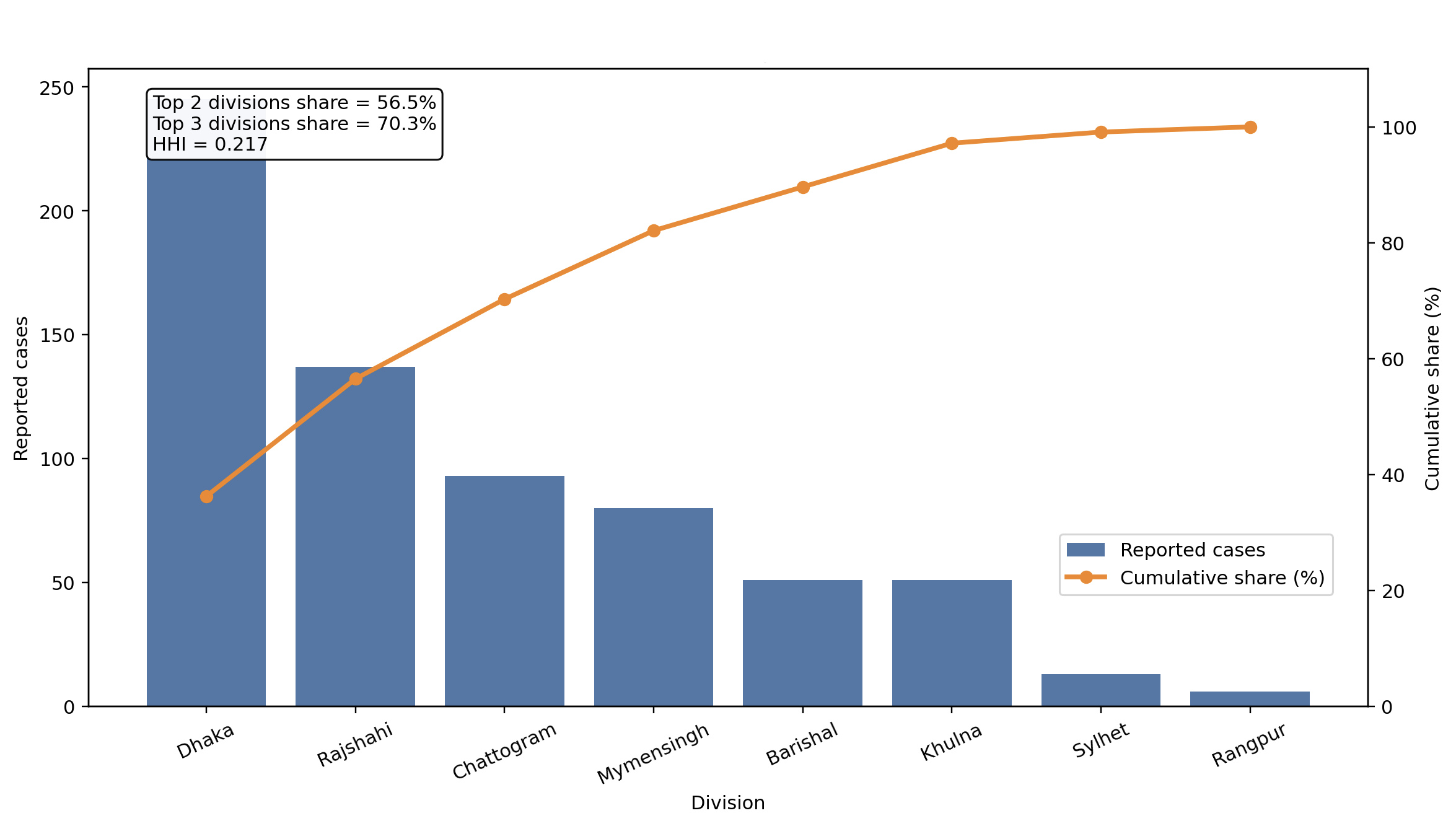}
\caption{Pareto analysis of division-wise measles reporting. Bars indicate division-level reported cases and the line shows the cumulative share of total burden in the reporting slice. The inset summarises the share held by the top two and top three divisions and the HHI value.}\label{fig:pareto_division}
\end{figure}

The outbreak burden was concentrated in a limited number of divisions rather than spread evenly across the country. Figure~\ref{fig:map_cases_deaths} shows broad national spread, but the 31 March 2026 division-level reporting slice showed Dhaka as the clear epicentre, with 245 reported cases, or 36.24\% of that slice. Rajshahi accounted for 137 cases and Chattogram for 93, meaning those three divisions together already dominated reported burden early in the surge. Figures~\ref{fig:division_slice} and \ref{fig:pareto_division} make the concentration explicit: the two highest-burden divisions accounted for 56.5\% of reported cases and the top three accounted for 70.3\%. The HHI was 0.217, indicating meaningful concentration rather than a diffuse national pattern. Operationally, that pattern justified a geographically targeted response rather than a flat allocation of outbreak resources \citep{reliefweb_bangladesh_measles_flash_2026,who_searo2026}

\subsection{Clinical and demographic profile}
\begin{figure}[H]
\centering
\includegraphics[width=0.88\textwidth]{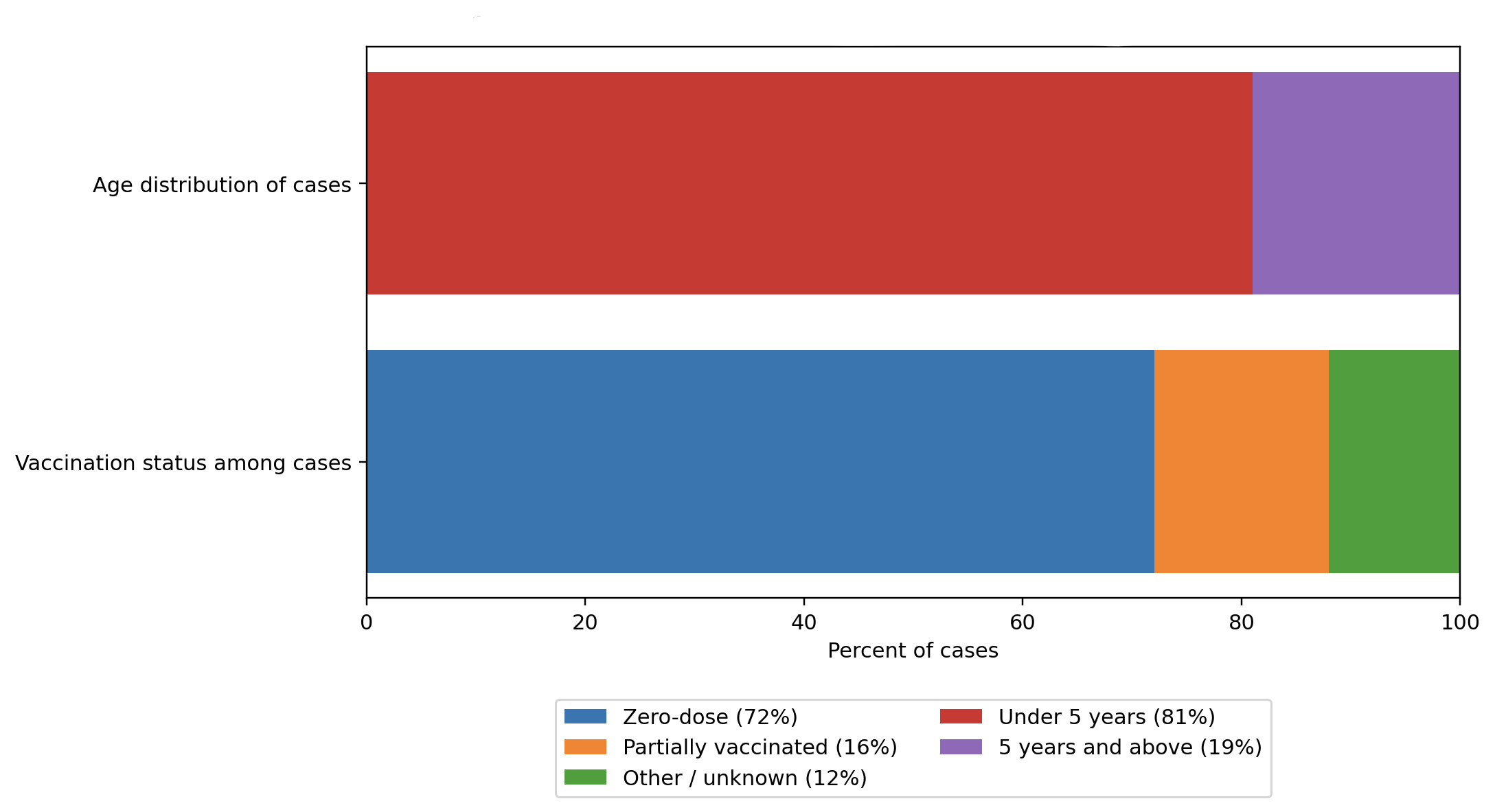}
\caption{Composition of the 2026 Bangladesh measles outbreak. The upper grouped bar shows the age distribution of reported cases and the lower grouped bar shows vaccination status among reported cases.}\label{fig:composition2026}
\end{figure}

The clinical profile of the outbreak was defined by very young children. UNICEF reported that 81\% of cases occurred in children younger than five years and 34\% in infants younger than nine months. Reported vaccination status reinforced the same message: 72\% of cases were zero-dose and 16\% were partially vaccinated. Figure~\ref{fig:composition2026} shows that the outbreak was concentrated in exactly the groups expected to be most vulnerable to transmission and severe outcomes when community protection fails. Hospital data added clinical weight to that pattern. At the Infectious Diseases Hospital in Dhaka, 685 suspected measles cases had been recorded, with 573 admissions in March 2026 compared with 69 in the same period in 2025. In Cox's Bazar, an eight-bed isolation unit was managing 30 to 40 children at a time. These figures show not just who was infected, but who was reaching hospitals under pressure \citep{UNICEF_SitRep1_2026,humanitarian}.

\subsection{Immunisation programme context and emergency response}
Programme indicators showed that Bangladesh entered 2026 with declining protection even before the outbreak accelerated. Partner reports described valid first-dose measles--rubella coverage falling from 88.6\% in 2019 to 86\% in 2024 and valid second-dose coverage from 89\% to 80.7\% over the same period. Those declines translated into large susceptible cohorts: roughly 10 million children on the first-dose gap estimate and about 20 million on the second-dose gap estimate. Figures~\ref{fig:trend_cases_coverage} %and \ref{fig:coverage_context} 
shows the longer trajectory. Response began on 5 April 2026 with an emergency campaign targeting more than 1.2 million children in 30 hot-spot upazilas across 18 districts; 75,442 children were vaccinated on the first day, equivalent to 13\% of the target population. The outbreak dose was extended to infants aged six months because transmission had already entered the pre-routine eligibility window. The response also had to compensate for another programme failure: only one of the two planned national vitamin A rounds had been completed in 2025, so the government prioritized remaining stocks for children with measles while additional procurement was sought \citep{WUENIC_BGD_2024,UNICEF_SitRep1_2026,WHO_BGD_Campaign_2026,humanitarian}.

%\begin{figure}[H]
%\centering
%\includegraphics[width=0.88\textwidth]{vc.png}
%\caption{Bangladesh measles vaccine coverage trends. The dashed horizontal line denotes the 95\% threshold commonly used as the benchmark for strong measles control.}\label{fig:coverage_context}
%\end{figure}

\section{Discussion}
\subsection{Why now? Interpreting the resurgence}
We interpret the 2026 resurgence as the predictable consequence of accumulated immunity gaps, not as evidence that the measles vaccine stopped working. The data support that reading from multiple directions. Figure~\ref{fig:historical_context} shows an outbreak far above Bangladesh's recent baseline. Figures~\ref{fig:division_slice} and \ref{fig:pareto_division} show concentration in a small number of divisions rather than uniform national spread. Figure~\ref{fig:composition2026} shows that cases clustered in zero-dose children, partially vaccinated children, and infants younger than nine months. Figures~\ref{fig:trend_cases_coverage} shows why that distribution was possible: second-dose coverage remained below the level required for durable interruption of measles transmission, while subnational inequity persisted behind national averages. That is the classic profile of a programme losing grip on population immunity. Similar dynamics have preceded major outbreaks elsewhere in the post-pandemic period, including in settings where national coverage looked respectable until local gaps widened enough for explosive transmission. Bangladesh's experience therefore fits a regional and global pattern. The failure was programmatic, not scientific \citep{LancetMicrobe_GlobalResurgence_2026,LancetSEA_Elimination_2025,WHO_measlesdeaths_2025,WHO_SEARO_EpiBull_2026}.

\subsection{Structural drivers and systems context}
Several structural failures turned vulnerability into crisis. First, routine services did not reach children equally. Urban slums, hard-to-reach districts, and refugee settlements have long required more intensive microplanning than national averages imply, and Cox's Bazar had already shown how quickly susceptibility can build in displaced populations. Second, the political transition of 2024 and 2025 appears to have interrupted routine programme continuity at precisely the moment when post-pandemic catch-up was needed most. Third, supply constraints mattered. Reports of depleted measles--rubella stocks and wider vaccine shortages meant that even motivated catch-up could not proceed at the necessary speed. Finally, the outbreak unfolded in a system where child nutrition protection had also weakened. Only one national vitamin A round was completed in 2025, increasing the likelihood that measles infection would translate into more severe disease. The lesson extends beyond Bangladesh. Nepal, Myanmar, Indonesia, and other SEAR countries face the same risk when uneven routine recovery, mobile populations, and crisis-sensitive delivery systems are allowed to coexist. In high-transmission diseases, those pressures do not remain local for long \citep{feldstein_rohingya_sero_2020,chin_rohingya_measles_model_2020,BusinessStandard_VaxShortage_2026,Guardian_Bangladesh_2026,humanitarian,WHO_SEARO_EpiBull_2026}.

\subsection{Research priorities and data requirements}
The next analytical step should be practical and immediate. Bangladesh needs upazila-level, line-listed surveillance data that capture age, vaccination history, nutritional status, complications, treatment pathway, and outcome. Without those data, it is impossible to estimate vaccine effectiveness credibly, identify the strongest predictors of severe disease, or separate transmission hotspots from referral hotspots. Viral genomic sequencing should also begin as soon as feasible to clarify whether multiple introductions occurred and to map transmission chains across districts and high-risk settings. These are not academic add-ons. They are the data required to run a better outbreak response in real time and to prevent the next one. Bangladesh already has methodological precedents for this work, including earlier vaccine-effectiveness studies and recent trial evidence supporting an outbreak dose at six months \citep{Akramuzzaman_BullWHO_2002,Ahmed_LancetChildAdolesc_2025}.

\subsection{Policy implications for Bangladesh and the region}
The policy implications are immediate. Bangladesh should continue the emergency measles--rubella campaign, but it should do so with sharper microgeographic targeting based on immunity-gap mapping rather than division-level totals alone. Vitamin A stocks must be replenished urgently, and paediatric surge capacity needs to expand in Dhaka, Rajshahi, and Chattogram where hospital pressure is already evident. Over the medium term, the country should restore the routine sessions disrupted since 2024, launch a national supplementary campaign before the next high-transmission period, and make outreach to urban slums and refugee populations a standing component of the immunisation strategy rather than an emergency exception. For the wider South-East Asia Region, Bangladesh should be treated as a sentinel event. IA2030 is already off track, and WUENIC 2025 shows that large cohorts of zero-dose and under-immunised children remain in lower-middle-income settings. WHO SEARO and national programmes need faster regional coordination before similar outbreaks become recurrent \citep{WHO_IA2030_2025,WUENIC_2025,WHO_BGD_Campaign_2026,WHO_SEARO_EpiBull_2026}.

%\section{Limitations}
%This analysis has four main limitations. First, it relied on aggregated rather than individual-level data, which means ecological fallacy remains a real risk and formal vaccine-effectiveness estimation was not possible. Second, the surveillance picture was shaped by reporting bias. Hospital-based and urban sources were more visible than rural community transmission, laboratory confirmation lagged behind clinical reporting, and community deaths were probably undercounted, all of which widened the gap between suspected and confirmed totals. Third, the evidence base reflected the constraints of an active outbreak. We used authoritative grey literature and carefully attributed media reports because those were often the fastest sources of named official figures and hospital observations, but real-time reporting inevitably contained discrepancies across dates, agencies, and case definitions. Vaccination histories reported in outbreak updates were also vulnerable to caregiver recall and incomplete records. Fourth, we had no genomic data. We could not examine viral lineage, trace introduction pathways, or assess whether more than one transmission chain was operating. Even with those limitations, the central epidemiological finding remained robust across sources: the 2026 Bangladesh outbreak was historically exceptional, geographically concentrated, and driven by accumulated immunity gaps rather than vaccine failure.

\section{Summary and Conclusion}
This analysis clarifies three things that were not previously assembled in one operational picture. First, the 2026 Bangladesh measles resurgence was not simply large; it was sharply concentrated in specific places and in specific children, with a small number of divisions carrying most of the early burden and with infants, zero-dose children, and partially vaccinated children bearing the heaviest risk. Second, the outbreak did not expose a failing vaccine. It exposed a strained programme. Declining second-dose coverage, missed children in urban and displaced settings, depleted stocks, and delayed catch-up created the conditions for explosive transmission. Third, measles severity in this outbreak cannot be separated from broader child-health erosion. Interrupted vitamin A delivery and health-system crowding widened the path from infection to hospitalisation and death.

These insights carry direct consequences for policy. Bangladesh cannot rely on national averages or one-off emergency messaging. It needs upazila-level immunity-gap mapping, sustained catch-up vaccination, restored routine delivery, and a permanent outreach strategy for urban slums, mobile populations, and refugee communities. It also needs better paediatric surge planning in high-burden hospitals and faster integration of surveillance, laboratory, and mortality review data into response operations. For South-East Asia more broadly, Bangladesh should be read as warning, not exception. Countries that allow post-pandemic susceptibility to accumulate will face the same reckoning. IA2030 will not be rescued by rhetoric or by national coverage headlines alone. Bangladesh, regional governments, WHO SEARO, UNICEF, Gavi, and partners must close the zero-dose gap now, before measles becomes once again a recurrent marker of preventable failure.

%To save lives and stop the spread of measles, health authorities and communities must take these urgent steps immediately.
%%	\subsection*{Immediate actions to reduce deaths}
%%Provide an actionable package:
%\begin{itemize}
%	\item Vaccination: outbreak-response MR vaccination and rapid catch-up for missed children.
%	\item Clinical care: early triage, pneumonia management, oxygen/respiratory support pathways, vitamin A treatment policy.
%	\item Infection control: isolation and prevention of nosocomial transmission.
%	\item Surveillance: strengthen case-based reporting, lab turnaround, and mortality audits.
%	\item Communication: counter misinformation; clarify vaccine effectiveness; promote timely care seeking.
%\end{itemize}

\section*{Declarations}
\textbf{Ethics approval and consent to participate:} This study was based exclusively on secondary analysis of aggregated, publicly available surveillance data published by the Directorate General of Health Services (Bangladesh), WHO, and UNICEF. No individual patient data were accessed. %In accordance with [Institution's] institutional policy on secondary data research, ethical review committee approval and individual participant consent were not required.

\medskip
\noindent\textbf{Consent for publication:} No consent is  required to publish this article. 

\medskip
\noindent\textbf{Availability of data and materials:} All data analysed in this study are derived from publicly accessible sources including DGHS Bangladesh daily surveillance bulletins (\url{dghs.gov.bd}), UNICEF Bangladesh measles outbreak situation reports (\url{unicef.org/bangladesh}), WHO SEARO epidemiological bulletins, and the WHO Global Immunization Data Portal (\url{immunizationdata.who.int}). No proprietary or restricted datasets were used.

\medskip
\noindent\textbf{Competing interests:} The authors declare no competing interests.

%\medskip
%\noindent\textbf{Funding:} [This research received no specific grant from any funding agency in the public, commercial, or not-for-profit sectors. / This work was supported by (grant name, funder, grant number).]

\medskip
\noindent\textbf{Authors' contributions (CRediT):} Conceptualisation: MK. Data curation: FK and AAS. Formal analysis:  FK, AAS,  NNKR and MK. Investigation: AAS and NNKR. Methodology: FK and MK. Visualisation: AAS. Writing -- original draft: FK, AAS and MK. Writing -- review \& editing: MK and NNKR. Supervision: MK. %Project administration: [Authors].

\medskip
\noindent\textbf{Acknowledgements:} The authors thank the Directorate General of Health Services, Bangladesh; UNICEF Bangladesh; WHO Bangladesh and WHO SEARO; and the Bangladesh Red Crescent Society for their publicly available situation reports and surveillance data that made this rapid analysis possible. The research was partially supported by the University Grants Commission (UGC), %the Ministry of Science and Technology,
Bangladesh.

\bibliographystyle{unsrtnat}
\bibliography{reference_cleaned_revised}

\end{document}